\documentclass[12pt]{article}
\usepackage{fullpage}
\usepackage{amssymb,amsmath,amsfonts,amsthm,latexsym,enumerate,url,cases}
\numberwithin{equation}{section}
\usepackage{hyperref}
\hypersetup{colorlinks=true,citecolor=blue,linkcolor=blue,urlcolor=blue}

\newtheorem{theorem}{Theorem}[section] %
\newtheorem{remark}[theorem]{Remark} %

\begin{document}
\title{On consecutive abundant numbers }

\author{  Yong-Gao Chen, Hui Lv\footnote{ygchen@njnu.edu.cn(Y.-G. Chen)}\\
\small  School of Mathematical Sciences and Institute of Mathematics, \\
\small  Nanjing Normal University,  Nanjing  210023,  P. R. China
}
\date{}
\maketitle \baselineskip 18pt \maketitle \baselineskip 18pt

{\bf Abstract.} A positive integer $n$ is called an abundant
number if $\sigma (n)\ge 2n$, where $\sigma (n)$ is the sum of all
positive divisors of $n$. Let $E(x)$ be the largest number of
consecutive abundant numbers not exceeding $x$. In 1935, P. Erd\H
os proved that there are two positive constants $c_1$ and $c_2$
such that $c_1\log\log\log x\le E(x)\le c_2\log\log\log x$. In
this paper, we resolve this old problem by proving that,
$E(x)/\log \log\log x$ tends to a limit as $x\to +\infty$, and the
limit value has an explicit form which is between $3$ and $4$.



\vskip 5mm

\section{Introduction}

Let $\sigma (n)$ is the sum of all positive divisors of $n$. A
positive integer $n$ is called an abundant number, a perfect
number and a deficient number if $\sigma (n)\ge 2n$, $=2n$ and
$<2n$, respectively. These numbers have brought extensive
research. For example, abundant numbers have been studied in
[1-15,17-19]. Let $E(x)$ be the largest number of consecutive
abundant numbers not exceeding $x$. In 1935, P. Erd\H os
\cite{Erdos1935b} proved that there are two positive constants
$c_1$ and $c_2$ such that $c_1\log\log\log x\le E(x)\le
c_2\log\log\log x$. P. Erd\H os paid much attention to abundant
numbers  all his life (see [11-15]).

In this paper, $p$ always denotes a prime and $(a, b)$ denotes the
greatest common divisor of two integers $a$ and $b$.

 The following result is
proved.

\begin{theorem} \label{thm1} We have $$\lim_{x\to \infty}\frac{E(x)}{\log \log\log x} =\left( \log \varpi \right)^{-1}, $$
where \begin{equation*}\varpi = \inf \left( \prod_{i=1}^M \max
\left\{ \frac{2(i, M)}{\sigma ((i, M))}, 1\right\} \right)^{\frac
1M}.\end{equation*} \end{theorem}

\begin{remark} From the proof of the main theorem, we also have
\begin{equation*}\varpi = \frac 1\beta \sup
\left( \prod_{i=1}^M \max\left\{ \frac{\sigma ((i, M))}{(i, M)},
2\right\}  \right)^{\frac 1M},\end{equation*} where
$$\beta = \prod_{p} \prod_{t=1}^\infty \left( 1+\frac 1p+\cdots +\frac
1{p^t} \right)^{\frac 1{p^t} \left( 1-\frac 1p\right) } .$$ We
have $\beta <1.56635$. For $M=4840909920000$, we have
$$ \varpi \ge \frac 1\beta \left( \prod_{i=1}^{M}
\max\left\{ \frac{\sigma ((i, M))}{(i, M)}, 2\right\}
\right)^{\frac 1{M}}>1.3267$$ and
$$\varpi \le \left( \prod_{i=1}^{M}
\max \left\{ \frac{2(i, M)}{\sigma ((i, M))}, 1\right\}
\right)^{\frac 1{M}}<1.3604.$$ It follows that
$$ 3.24<\left( \log \varpi \right)^{-1} < 3.54.$$
It is not the aim of this paper to find a good numerical result.
\end{remark}

\medspace

Now we give an outline of the proof of Theorem \ref{thm1}. One of
the main ideas in this paper is to introduce the following sequence:
$$\{ (n, M)\}_{n=1}^\infty$$ for any given positive integer $M$. This
is a  sequence of period length $M$. If $M$ is divisible by all
``small" prime powers, then two sequences
$$\left\{ \frac{\sigma (n)}{n} \right\}_{n=1}^\infty
\text{  and  } \left\{ \frac{\sigma ((n, M))}{(n, M)} \right\}_{n=1}^\infty
$$
are ``similar" in a sense, but the late one is a sequence of
period length $M$. I believe that the sequence $\{ (n,
M)\}_{n=1}^\infty$ will play an important role  in the future
research.

 In Section
\ref{secupper}, we give an upper bound of $E(x)$: $$E(x)\le
(\rho_1 (M)+o(1)) \log \log\log x$$ for a function $\rho_1 (M)$
and any given positive integer $M$. In Section \ref{seclower}, we
give a lower bound of $E(x)$: $$E(x)\ge (\rho_2 (M)+o(1)) \log
\log\log x$$ for a function $\rho_2 (M)$ and any given positive
integer $M$. In Section \ref{secrelation}, it is proved that
$$\rho_1 (M)=\rho_2 (M) +o(1)$$ for infinitely many positive
integers $M$. In Section \ref{uniform}, let $\rho_1 =\inf \rho_1
(M)$ and $\rho_2 =\sup \rho_2 (M)$. We finish the proof by proving
that $\rho_1=\rho_2=(\log \varpi )^{-1}$.

Let $\alpha $ be a positive real number. A positive integer $n$ is
called an $\alpha$-abundant number if $\sigma (n)\ge  \alpha n$.
Let $E_\alpha (x)$ be the largest number of consecutive
$\alpha$-abundant numbers not exceeding $x$.

The method in this paper can be used to prove the following
analogous result.

\begin{theorem} \label{thm2} We have $$\lim_{x\to \infty}\frac{E_\alpha (x)}{\log \log\log x}
=\left( \log \varpi_\alpha \right)^{-1}, $$ where
\begin{equation*}\varpi_\alpha = \inf \left( \prod_{i=1}^M \max
\left\{ \frac{ (i, M)\alpha}{\sigma ((i, M))}, 1\right\}
\right)^{\frac 1M}.\end{equation*} \end{theorem}

\section{The upper bound of $E(x)$}\label{secupper}

In this section we prove that, for any given positive integer $M$,
$$E(x)\le (\rho_1 (M) +o(1)) \log \log\log x $$ for all sufficient large
$x$, where
$$\rho_1 (M)=(\log \delta_M -\log \beta )^{-1},$$
$$\beta = \prod_{p} \prod_{t=1}^\infty \left( 1+\frac 1p+\cdots +\frac
1{p^t} \right)^{\frac 1{p^t} \left( 1-\frac 1p\right) }$$ and
$$\delta_M=\left( \prod_{i=1}^{M}\max\left\{ \frac{\sigma ((i,
M))}{(i, M)}, 2\right\}\right)^{\frac 1M}.$$ It is clear that
$$\beta = \prod_{p} \prod_{t=1}^\infty \left( \frac{\sigma (p^t)}{p^t} \right)^{\frac 1{p^t} \left( 1-\frac 1p\right) }.$$

Let $m, m+1, \dots , m+k-1$ be consecutive abundant numbers not
exceeding $x$. Let $M$ be a positive integer. For any prime $p$ and
any positive integer $t$, let $s_{p^t}$ be the number of integers in
$m, m+1, \dots , m+k-1$ which are divisible by $p^t$. Then
\begin{equation}\label{eq0}s_{p^t}=\frac {k}{p^t} +R_{p^t},\quad |R_{p^t}|\le
1.\end{equation} If $p^t>m+k-1$, then $s_{p^t}=0$.

For any two positive integers $a,b$ with $a\mid b$, we have
$$\frac{\sigma (a)}{a}\le \frac{\sigma (b)}{b}.$$
Thus
$$\frac{\sigma ((m+i, M))}{(m+i, M)}\le \frac{\sigma (m+i)}{m+i},
\quad i=0, 1, \dots , k-1.$$ Since $m+i$ $(1\le i\le k-1)$ are
abundant numbers, we have
$$2\le \frac{\sigma (m+i)}{m+i},
\quad i=0, 1, \dots , k-1.$$ It follows that
$$\prod_{i=0}^{k-1}
\max\left\{ \frac{\sigma ((m+i, M))}{(m+i, M)}, 2 \right\} \le
\prod_{i=0}^{k-1} \frac{\sigma (m+i)}{m+i}.$$ By the definition of
$s_{p^t}$, we have
$$\prod_{i=0}^{k-1} \frac{\sigma (m+i)}{m+i} \\
=\prod_{p\mid m(m+1)\cdots (m+k-1)} \prod_{t=1}^\infty \left(
\frac{\sigma (p^t)}{p^t}\right)^{s_{p^t}-s_{p^{t+1}}}.$$ If $p^t\le
k$, then $p\mid m(m+1)\cdots (m+k-1)$. Now we split the left product
into two parts according to $p^t\le k$ and $p^t>k$:
\begin{equation}\label{eqx} \prod_{p} \prod_{t=1, p^t\le k}^\infty \left(
\frac{\sigma (p^t)}{p^t}\right)^{s_{p^t}-s_{p^{t+1}}} \end{equation}
and  $$\prod_{ p\mid m(m+1)\cdots (m+k-1)} \prod_{t=1, p^t>k}^\infty
\left( \frac{\sigma (p^t)}{p^t}\right)^{s_{p^t}-s_{p^{t+1}}}.
$$
 It is clear that
$$\frac{\sigma (p^t)}{p^t}=1+\frac 1p+\cdots +\frac
1{p^t}<1+\frac 1{p-1}.$$ Let $T=T_p$ be the integer with
$p^{T_p}\le k<p^{T_p+1}$. Since there is at most one integer in
$m, m+1, \dots , m+k-1$  which is divisible by $p^{T_p+1}$, it
follows that $s_{p^{T_p+1}}\le 1$. Thus
\begin{eqnarray*}&&\prod_{p\mid m(m+1)\cdots (m+k-1)} \prod_{t=1, p^t>k}^\infty \left( \frac{\sigma
(p^t)}{p^t}\right)^{s_{p^t}-s_{p^{t+1}}} \\
&\le &   \prod_{ p\mid m(m+1)\cdots (m+k-1)} \prod_{t=1,
p^t>k}^\infty \left( 1+\frac
1{p-1}\right)^{s_{p^t}-s_{p^{t+1}}}\\
&=&  \prod_{ p\mid m(m+1)\cdots (m+k-1)} \left( 1+\frac
1{p-1}\right)^{s_{p^{T_p+1}}} \\
 &\le &  \prod_{ p\mid m(m+1)\cdots
(m+k-1)} \left( 1+\frac 1{p-1}\right)\\
&=&\frac{m(m+1)\cdots
(m+k-1)}{\phi (m(m+1)\cdots (m+k-1))}.
\end{eqnarray*}
Noting that (see \cite[Theorem 328]{Hardy1979}),
$$\phi (n)\gg \frac{n}{\log\log n},$$
where $\gg$ is the Vinogradov symbol, we have
\begin{eqnarray*} &&\prod_{p\mid m(m+1)\cdots (m+k-1)} \prod_{t=1, p^t>k}^\infty \left( \frac{\sigma
(p^t)}{p^t}\right)^{s_{p^t}-s_{p^{t+1}}}\\
&\le & \frac{m(m+1)\cdots (m+k-1)}{\phi (m(m+1)\cdots (m+k-1))}\\
&\ll & \log\log (m(m+1)\cdots (m+k-1))\nonumber \\
&\ll & \log\log (x^k)\\
&\ll & \log k+\log\log x\\
&\ll & (\log k) \log\log x,
\end{eqnarray*}
where $\ll$ is also the Vinogradov symbol.

Now we deal with the product \eqref{eqx}.

By \eqref{eq0}, we have
\begin{eqnarray*}s_{p^t}-s_{p^{t+1}}&=&\frac {1}{p^t} \left( 1-\frac 1p\right)
k +R_{p^t}-R_{p^{t+1}}\\
&\le& \frac {1}{p^t} \left( 1-\frac 1p\right) k +2.\end{eqnarray*}
Noting that $T_p=0$ for $p>k$ and $T_p\le 2\log k$ for any prime
$p$, we have
\begin{eqnarray*}&&\prod_{p} \prod_{t=1, p^t\le k}^\infty \left(
\frac{\sigma (p^t)}{p^t}\right)^{s_{p^t}-s_{p^{t+1}}}\\
&\le & \prod_{p\le k} \prod_{t=1}^{T_p} \left( \frac{\sigma
(p^t)}{p^t}\right)^{\frac {1}{p^t} \left( 1-\frac 1p\right) k +2}\\
&\le & \prod_{p} \prod_{t=1}^\infty \left( \frac{\sigma
(p^t)}{p^t}\right)^{\frac {1}{p^t} \left( 1-\frac 1p\right) k} \cdot
\prod_{p\le k } \prod_{t=1}^{T_p} \left( \frac{\sigma
(p^t)}{p^t}\right)^{2}\\
&\le & \prod_{p} \prod_{t=1}^\infty \left( \frac{\sigma
(p^t)}{p^t}\right)^{\frac {1}{p^t} \left( 1-\frac 1p\right) k} \cdot
\prod_{p\le k }  \left( 1+\frac 1{p-1} \right)^{4\log k}\\
&=& \beta^k  \prod_{p\le k }\left( 1+\frac 1{p-1} \right)^{4\log k}\\
&\le & \beta^k  (c\log k)^{4\log k},
\end{eqnarray*}
where the last inequality is due to Mertens' theorem (see
\cite[Theorem 429]{Hardy1979}) and $c$ is a positive constant.
Hence
\begin{eqnarray*}\prod_{i=0}^{k-1} \max\left\{ \frac{\sigma ((m+i, M))}{(m+i, M)}, 2
\right\} \ll  \beta^{k} (c\log k)^{4\log k+1} \log\log x.
\end{eqnarray*}
Recall that
$$\delta_M=\left( \prod_{i=1}^{M}\max\left\{ \frac{\sigma ((i,
M))}{(i, M)}, 2\right\}\right)^{\frac 1M}.$$ Since
\begin{eqnarray*}&&\prod_{i=0}^{k-1} \max\left\{ \frac{\sigma ((m+i, M))}{(m+i, M)}, 2
\right\}\\
&\ge &  \left( \prod_{i=1}^{M}\max\left\{ \frac{\sigma ((i,
M))}{(i, M)}, 2\right\}\right)^{k/M-1}\\
&\gg_M& \delta_M^k ,\end{eqnarray*} it follows that
$$\delta_M^k \ll_M  \beta^{k} (c\log k)^{4\log k+1} \log\log x.$$
Noting $\beta <2\le \delta_M$, we have
\begin{eqnarray*}k&\le & (\log \delta_M -\log \beta +o(1))^{-1} \log \log\log x\\
&=&(\rho_1 (M) +o(1))\log \log\log x.\end{eqnarray*}

\section{The lower bound of $E(x)$} \label{seclower}

In this section we prove that, for any given positive integer $M$,
$$E(x)\ge ( \rho_2 (M) +o(1) ) \log \log\log
x$$ for all sufficient large $x$, where
$$\rho_2 (M)=\left( \log 2 -\log \tau_M \right)^{-1},$$
$$\tau_M =\left( \prod_{i=1}^{M}\min\left\{ \frac{\sigma ((i, M))}{(i,
M)}, 2\right\}\right)^{1/M}.$$

Let $M$ be any given positive integer and let $q_1, q_2, \dots $
be all primes in increasing order which are greater than $M$. Let
$$A=\prod_{M<p<\frac 12\log x} p.$$
 By Mertens' theorem (see \cite[Theorem 429]{Hardy1979}), we have
\begin{equation}\label{eqy}\frac{\sigma (A)}{A} =\prod_{M<p<\frac 12\log x} \left( 1+ \frac 1p \right) \ge
c_M \log\log x,\end{equation} where $c_M$ is a positive constant
depending only on $M$.

 Let
$j_0=0$. For any integer $l\ge 1$, let
 $$a_l=(l, M) q_{j_{l-1}+1} \cdots q_{j_l} := (l, M) b_l,$$ where $j_l$ is the least integer
with $j_l\ge j_{l-1}+1$ and $\sigma (a_l)\ge 2 a_l$.

Since $((l,M), b_l)=1$, it follows that
$$\frac{\sigma (a_l)}{a_l} = \frac{\sigma (b_l)}{b_l} \frac{\sigma ((l, M))}{(l,
M)}$$
It is clear that
\begin{eqnarray*}\frac{\sigma (a_l)}{a_l}&<&\max\left \{ 2, \frac{\sigma ((l, M))}{(l, M)} \right\}   \left( 1+
\frac{1}{q_{j_l}}\right)\\
&\le & \max\left \{ 2, \frac{\sigma ((l, M))}{(l, M)} \right\}
\left( 1+ \frac{1}{l}\right) .\end{eqnarray*} Let $k$ be the integer
with
$$b_1b_2\cdots b_k\le A< b_1b_2\cdots b_{k+1}.$$
Then
\begin{eqnarray*}\frac{\sigma (A)}{A}&<& \frac{\sigma (b_1)}{b_1} \cdots \frac{\sigma
(b_{k+1})}{b_{k+1}}\\
&=& \frac{\sigma (a_1)}{a_1} \cdots \frac{\sigma
(a_{k+1})}{a_{k+1}}\frac{(1, M)}{\sigma ((1, M))} \cdots
\frac{(k+1, M)}{\sigma ((k+1, M))} \\
&\le & \prod_{i=1}^{k+1} \max\left \{ 2, \frac{\sigma ((i,
M))}{(i, M)}
\right\} \cdot  \prod_{i=1}^{k+1} \left( 1+ \frac{1}{i}\right)\cdot \prod_{i=1}^{k+1} \frac{(i, M)}{\sigma ((i, M))}\\
&=& 2^{k+1} (k+2) \prod_{i=1}^{k+1}\max\left\{ \frac{(i,
M)}{\sigma
((i, M))}, \frac 12\right\}\\
&\le & 2^{k+1} (k+2) \left( \prod_{i=1}^{M}\max\left\{ \frac{(i,
M)}{\sigma ((i, M))}, \frac 12\right\}\right)^{(k+1)/M-1}.
\end{eqnarray*}
Recall that
$$\tau_M =\left( \prod_{i=1}^{M}\min\left\{ \frac{\sigma ((i, M))}{(i,
M)}, 2\right\}\right)^{1/M},$$ we have
\begin{eqnarray*}\frac{\sigma (A)}{A}&\le& 2^{k+1} (k+2) \left( \prod_{i=1}^{M}\max\left\{ \frac{(i,
M)}{\sigma ((i, M))}, \frac 12\right\}\right)^{(k+1)/M-1}\\
&\ll_M& 2^{k+1} k \tau_M^{-k}.\end{eqnarray*} It follows from
\eqref{eqy} that
$$c_M \log\log x \le \frac{\sigma (A)}{A} \ll_M
2^{k+1} k \tau_M^{-k}.$$ Noting that $\tau_M <2$, we have
\begin{eqnarray*}k&\ge& \left( \log 2 -\log \tau_M +o(1) \right)^{-1} \log \log\log
x\\
&=& (\rho_2 (M) +o(1))\log \log\log x.\end{eqnarray*}

Although we can prove that $k\ll \log\log\log x$, in order to
avoid unnecessary arguments, we prefer to write
\begin{equation}\label{eqw}k'=\min\{ k,\ 2 \lfloor \rho_2 (M) \log\log \log x\rfloor
\} ,\end{equation} where $\lfloor z\rfloor$ denotes the integral
part of real number $z$.

Now we prove that there are $k'$ consecutive abundant numbers not
exceeding $x$. It follows that
$$E(x)\ge (\rho_2 (M) +o(1))\log \log\log x.$$

By the Chinese remainder theorem (see \cite[Theorem
121]{Hardy1979}), there exists a positive integer $ m\le M
b_1\cdots b_{k'}$ such that $m\equiv 0\pmod{M}$ and
$$m+i\equiv 0\pmod{b_i}, \quad i=1,2, \dots, k'.$$

Now we prove that $m+1, m+2, \dots , m+k'$ are consecutive
 abundant numbers which do not exceed $x$.

By the prime number theorem (see \cite[Theorem 6]{Hardy1979}), we
have
$$\log (k'M A) <k'M+\log A<k'M+\frac 23 \log x<\log x$$
for all sufficiently large $x$. It follows that $$m+k'\le mk'\le
k'Mb_1\cdots b_{k'}\le k'MA<x$$ for all sufficiently large $x$.

For $1\le i\le k'$, by $m\equiv 0\pmod{M}$ we have $(i, M)\mid m$
and then $(i, M)\mid m+i$. Noting that $b_i\mid m+i$ and $((i, M),
b_i)=1$, we have $(i, M)b_i \mid m+i$. It follows that $a_i \mid
m+i$. Since $a_i$ $(1\le i\le k')$ are abundant numbers, it
follows that $$m+1, m+2, \dots , m+k'$$ are consecutive
 abundant numbers.

\section{$\rho_1 (M)  =\rho_2 (M) +o(1)$ for infinitely many
$M$}\label{secrelation}

Let $U$ be a large integer and let
$$M_U =\prod_{p<U} p^U.$$

In this section we  prove that
\begin{equation}\label{eqm}\rho_1 (M_U)  =\rho_2 (M_U) +o(1).\end{equation}
Recall that
$$\rho_1 (M)=(\log \delta_M -\log \beta )^{-1}$$ and
$$\rho_2 (M)=\left( \log 2 -\log \tau_M \right)^{-1},$$
it is enough to prove that
\begin{equation}\label{eqz}\delta_{M_U} \tau_{M_U} =2\beta +o_U(1).\end{equation}

Let $M=M_U$. For any prime $p<U$ and any positive integer $t\le
U$, let $v_{p^t}$ be the number of integers in $1, 2, \dots , M$
which are divisible by $p^t$. Then
\begin{equation*}v_{p^t}=\frac {M}{p^t}.\end{equation*}
Let $v_{p^{U+1}}=0$ for any prime $p<U$.
 Thus
\begin{eqnarray*}\delta_{M_U} \tau_{M_U}&=& \left( \prod_{i=1}^{M}\max\left\{ \frac{\sigma ((i,
M))}{(i, M)}, 2\right\}\right)^{\frac 1M} \left(
\prod_{i=1}^{M}\min\left\{ \frac{\sigma ((i, M))}{(i, M)},
2\right\}\right)^{\frac 1M}\\
&=& \left( \prod_{i=1}^{M}\max\left\{ \frac{\sigma ((i, M))}{(i,
M)}, 2\right\}  \min\left\{ \frac{\sigma ((i, M))}{(i, M)},
2\right\} \right)^{\frac 1M}\\
&=& \left( \prod_{i=1}^{M}\left( 2\frac{\sigma ((i, M))}{(i,
M)}\right) \right)^{\frac 1M}\\
&=& 2\left( \prod_{i=1}^{M}\frac{\sigma ((i, M))}{(i,
M)} \right)^{\frac 1M}\\
&=& 2 \left(\prod_{p<U} \prod_{t=1}^U \left( \frac{\sigma
(p^t)}{p^t}\right)^{v_{p^t}-v_{p^{t+1}}}\right)^{\frac 1M}\\
&=&2 \prod_{p<U} \prod_{t=1}^U \left( 1+\frac 1p+\cdots +\frac
1{p^t} \right)^{\frac 1{p^t} \left( 1-\frac 1p\right) }\\
&& \cdot \prod_{p<U} \left( 1+\frac 1p+\cdots +\frac 1{p^{U+1}}
\right)^{\frac 1{p^{U+1}}
} \\
&=& 2\beta +o_U (1),
\end{eqnarray*}
where
\begin{eqnarray*}&& \log \prod_{p<U} \left( 1+\frac 1p+\cdots +\frac 1{p^{U+1}}
\right)^{\frac 1{p^{U+1}} }\\
&<& \sum_{p<U} \frac 1{p^{U+1}} \log \left(1+\frac 1{p-1}\right)\\
&<& \sum_{p<U} \frac 1{p^{U+1}(p-1)}\\
&<& \frac 1{p^{U}} \sum_{p} \frac 1{p(p-1)}\\
&=& o_U(1)
\end{eqnarray*}
and then
$$\prod_{p<U} \left( 1+\frac 1p+\cdots +\frac 1{p^{U+1}}
\right)^{\frac 1{p^{U+1}} }=1+o_U(1).$$ Thus we have proved that
\eqref{eqz} holds and so does \eqref{eqm}.

\section{Proof of Theorem \ref{thm1}}\label{uniform}

We have proved that, for any given positive integer $M$,
$$E(x)\le (\rho_1 (M) +o(1)) \log \log\log x $$ and
$$E(x)\ge  ( \rho_2 (M) +o(1) ) \log \log\log
x$$ for all sufficient large $x$. In order to obtain  the optimal
upper bound and the optimal lower bound of $E(x)$, we should
choose $M_1$ and $M_2$ such that $\rho_1 (M_1)$ is as small as
possible and $\rho_2 (M_2)$  is as large as possible. Let $$\rho_1
=\inf \rho_1 (M),\quad \rho_2 =\sup \rho_2 (M).
$$
Then
$$\rho_2+o(1)\le \frac{E(x)}{ \log \log\log x} \le  \rho_1+o(1).$$
So $\rho_2\le \rho_1$.

Now we prove that $\rho_2\ge \rho_1$.

Let $U$ be a large integer and $M_U$ be as in the previous section.
Then
$$\rho_1\le \rho_1 (M_U),\quad \rho_2\ge \rho_2 (M_U).$$
Since
$$\rho_1 (M_U)  =\rho_2 (M_U) +o(1),$$
it follows that
 $$\rho_1\le \rho_1 (M_U)=\rho_2 (M_U) +o(1)\le \rho_2+o(1).$$ This implies that $\rho_1\le \rho_2$.
Therefore, $\rho_2= \rho_1$ and then
$$\frac{E(x)}{ \log \log\log x}=\rho_2+o(1)=( \log \varpi )^{-1}+o(1).$$
This completes the proof of our main theorem.

\noindent\textbf{Acknowledgments.}  The first author is supported
by the National Natural Science Foundation of China, Grant No.
11371195 and a project funded by the Priority Academic Program
Development of Jiangsu Higher Education Institutions.


\begin{thebibliography}{99}


\bibitem{AlaogluErdos1944} L. Alaoglu and P. Erd\H os,  On highly composite and similar numbers,
Trans. Amer. Math. Soc. 56, (1944), 448--469.

\bibitem{Avidon1996} M. R. Avidon, On the distribution of primitive abundant
numbers,
Acta Arith. 77 (1996), 195--205.

\bibitem{Behrend1932}F. Behrend, \" Uber numeri abundantes, Preuss. Akad. Wiss.
Sitzungsber 21/23 (1932), 322--328.

\bibitem{Behrend1933}F. Behrend, \" Uber numeri abundantes, II, Preuss. Akad. Wiss.
Sitzungsber 6 (1933), 280--293.

\bibitem{Briggs2006}K. Briggs, Abundant numbers and the Riemann
hypothesis, Experiment. Math. 15 (2006),  251--256.

\bibitem{Chowla1934} S. Chowla, On abundant numbers, J. Indian Math. Soc. 1 (1934),
41--44.

\bibitem{Cohen1984} G. L. Cohen,  Primitive $\alpha $-abundant numbers, Math. Comp. 43 (1984),
 263--270.

\bibitem{Davenport1933} H. Davenport, จน\" uber numeri abundantes, S.-Ber. Preu. Akad.
Wiss., math.-nat. Kl. (1933), 830--837.


\bibitem{Dickson1913-413} L. E. Dickson,  Finiteness of the odd perfect and primitive
abundant numbers with $n$ distinct prime factors,  Amer. J. Math. 35
(1913),  413--422.


\bibitem{Dickson1913-423} L. E. Dickson,  Even abundant numbers,  Amer. J. Math. 35 (1913), 423--426.



\bibitem{Erdos1934} P. Erd\H os,  On the density of the abundant numbers,  J. London
Math. Soc. 9 (1934), 278--282.


\bibitem{Erdos1935a} P. Erd\H os, On primitive abundant numbers, J. London Math. Soc. 10
(1935), 49--58.


\bibitem{Erdos1935b} P. Erd\H os, Note on consecutive abundant
numbers, J. London Math. Soc. 10 (1935), 128--131.



\bibitem{Erdos1935} P. Erd\H os,  Remarks on number theory. I. On primitive $\alpha $-abundant
numbers,  Acta Arith. 5 (1958), 25--33.


\bibitem{Erdos1974} P. Erd\H os, On abundant-like numbers, Canad. Math. Bull. 17 (1974),
 599--602.

\bibitem{Hardy1979} G. H. Hardy and E. M. Wright, An Introduction to
the theory of numbers, Oxford Univ. Press 1979.

\bibitem{Kobayashi2014}M. Kobayashi, A new series for the density of abundant
numbers, Int. J. Number Theory 10 (2014),  73--84

\bibitem{Pollack2014} M. Kobayashi and P. Pollack, The error term in the count of abundant numbers, Mathematika 60 (2014),
 43--65.

\bibitem{Shapiro1968} H. N. Shapiro,  On primitive abundant numbers, Comm. Pure Appl.
Math. 21 (1968), 111--118.





\end{thebibliography}
\end{document}